\documentclass{amsart}
\usepackage{amsmath}%
\usepackage{amsfonts}%
\usepackage{amssymb}%
\usepackage{graphicx}
\theoremstyle{plain}

\numberwithin{equation}{section}
\begin{document}
\title[On Effective Lower Bounds for Points with  Dense Orbit]{  Canonical Heights and Monomial Maps: On Effective Lower Bounds for Points with  Dense Orbit}
\author{Jorge Mello}
\address[]
{Universidade Federal do Rio de Janeiro, Instituto de Matem\'{a}tica. mailing adress:\newline Rua Aurora 57/101, Penha Circular, 21020-380 Rio de Janeiro, RJ, Brasil. %
 } \email[]{jmelloguitar@gmail.com}%
\urladdr{https://sites.google.com/site/algebraufrj/}

\thanks{}
\date{Januray 8, 2019.}
\dedicatory{}

\begin{abstract}
  We prove, for the canonical height defined by Silverman [15] on monomial maps, the existence of effective lower bounds for heights of points with Zariski dense orbit, for cases with endomorphisms induced by matrices with real Jordan form.
\end{abstract}
\maketitle

\section{Introduction}

Height functions are used to study the arithmetic complexity of rational points on defined algebraic varieties. They have being important in the development of Diophantine geometry, and more recently, of arithmetic dynamics. The theory of canonical heights associated with morphisms $\phi : \mathbb{P}^N \rightarrow \mathbb{P}^N$ is well known [4].

For a system $(X/K, f_1,...,f_k, L)$ with $k$ self-morphisms on a smooth algebraic variety $X$ over a number field $K$, and $L$ a divisor satisfying a linear equivalence $\otimes^k_{i=1}f^*_i(L) \sim L^{\otimes d}$ for $d>k$, there is also a known theory of canonical heights developed by Kawaguchi in [8].

For rational maps, it is harder to construct and to study canonical height functions. In [7], Kawaguchi studied and understood the case of regular affine automorphisms. Recently [15], Silverman constructed and studied the theory of canonical height functions for dominant rational maps, and in particular, their behaviour for the case of monomial maps.

Namely, for a dominnt rational map $\phi : \mathbb{P}^N \dashrightarrow \mathbb{P}^N$, the first dynamical degree is defined as
\begin{center}
$\delta_{\phi}= \lim_{n \rightarrow \infty} \deg( \phi^n)^{\frac{1}{n}}$
\end{center}

Assuming that $\phi$ is defined over $\bar{\mathbb{Q}}$, there is a strong conjecture claiming that the infimum 
\begin{center}
$l_{\phi}:= \inf \{ l \geq 0 : \sup_{n \geq 1} \dfrac{\deg( \phi^n)}{n^l \delta_{\phi}^n} < \infty \}$
\end{center} exists and it is an integer such that $0 \leq l_{\phi} \leq N$ [15]. The set $\mathbb{P}^N(\bar{\mathbb{Q}})_{\phi}$ is defined as the subset of $\mathbb{P}^N(\bar{\mathbb{Q}})$ whose points $P$ have forward $\phi-$orbit $\mathcal{O}_{\phi}(P)$ well defined. If the conjecture is true, then, for each $P \in \mathbb{P}^N(\bar{\mathbb{Q}})_{\phi}$, Silverman defined the canonical height of $P$ with respect to $\phi$ by
\begin{center}
$ \hat{h}_{\phi}(P)=\lim \sup_{n \rightarrow \infty} \dfrac{1}{n^{l_{\phi}} \delta_{\phi}^n} h(\phi^n(P))$,
\end{center} where $h$ denotes the usual Weil height (see [6]).

For general dominant rational maps, the conjecture is open. However, for monomial maps, the conjecture is true, and the degree sequence is understood. Given a matrix $A=(a_{ij})$ with integer coordinates, one associates a self map $\phi_A$ on the algebraic multiplicative group, the algebraic torus $\mathbb{G}^N_m(\bar{\mathbb{Q}})$, by the following:
\begin{center}
$\phi_A(x_1,...,x_N)=({x}_{1}^{a_{11}}{x}_{2}^{a_{12}}...{x}_{N}^{a_{1N}},..., {x}_{1}^{a_{N1}}{x}_{2}^{a_{N2}}...{x}_{N}^{a_{NN}}).$
\end{center}
The group endomorphism $\phi_A: \mathbb{G}^N_m(\bar{\mathbb{Q}}) \rightarrow \mathbb{G}^N_m(\bar{\mathbb{Q}})$ is called the monomial map associated to $A$. It is surjective if $\det A \neq 0$, and it extends to a rational self map on $\mathbb{P}^N(\bar{\mathbb{Q}})$ denoted by $\phi_A$.
The dynamical degree of $\phi_A$ is equal to the spectral radius of $A$ (see [15, Proposition 21]) and the value $l_{\phi_A}$ is an integer satisfying $0 \leq l_{\phi_A} \leq N-1$, determined by its Jordan form (see 15, Theorem 24). The canonical height function for monomial maps is studied in [15]. In this paper, Silverman proved some results for monomial maps on an algebraic torus, as follows: \newline 

{\bf Theorem 1.1:} [Silverman, 15] \textit{Let $\phi_A$ be a monomial map with $\delta_{\phi_A} >1$.}

 (a) \textit{Let $P \in \mathbb{G}^N_m(\bar{\mathbb{Q}})$ be a point with $\hat{h}_{\phi_A}(P)=0.$ Then $\mathcal{O}_{\phi_A}(P)$ is contained in a proper non Zariski dense algebraic subgroup of $\mathbb{G}^N_m$.}

(b) \textit{If the characteristic polynomial of the matrix $A$ is irreducible over $\mathbb{Q}$, then}
\begin{center}
\textit{$\hat{h}_{\phi_A}(P)=0 \iff \mathcal{O}_{\phi_A}(P)$ is finite.}
\end{center}

According to all this, it follows that a point with Zariski dense orbit must have positive canonical height. A subsequent natural question then, made by Silverman [Remark 30, 15], is to ask if one can find explicitely a positive lower bound, effectively computable, for the canonical height of a point that has Zariski dense orbit, where the bound would depend only on the naive height of the point and on the matrix which induces the endomorphism. Silverman's suggestion is to try to use some effective form of a famous theorem of Alan Baker for such, considering that a conceptual form of it was used to prove the Theorem above. 

A goal for this paper is to answer Silverman's question, confirming his predictions for maps induced by matrices with real eigenvalues and points with algebraic integer coodinates, and also to obtain similar results for more general cases as those with complex eigenvalues, as far as one assumes that some information about certain limits of powers of the eigenvalues is known . More precisely, we prove the following:
\newline \newline {\bf Theorem A:}
\textit{Let $\phi:  {\mathbb{G}}_{m}^N \rightarrow  {\mathbb{G}}_{m}^N$ be a monomial map induced by $A$ with real eigenvalues, $l(A) \geq 1$ and $\delta_{\phi}>1$, and let $P \in  {\mathbb{G}}_{m}^N(\bar{\mathbb{Q}})$ be a point with Zariski dense orbit $\mathcal{O}_{\phi}(P)$. Let $J(A)$ be the matrix whose columns form a basis in which $A$ has its Jordan form, and let $\mathbb{Q}(J(A),P)$ be a field of definition of $P$ and $J(A)$ which is Galois over $\mathbb{Q}.$ Suppose that the coordinates of $P$ are algebraic integers in $\mathbb{Q}(J(A),P)$. Then there is an effective computable positive constant $C$ depending on $J(A)$, on $\mathbb{Q}(J(A),P)$, and on $h(P)$, such that $ \hat{h}_{\phi_{A}}(P)>C$}.
\newline \newline For this, we use an algorithm to compute Jordan normal forms from [14], and an effective form of Baker's theorem due to Philippon and Waldschmidt. As Silverman predicted, we see that the bigger the naive height of the point is, the smaller the calculated contants wil be. As a particular case, we obtain also effectivity for computing the lower bound of points with infinite orbit, when the matrix has irreducible characteristic polynomial.
\newline \newline 
{\bf Corollary B:}
\textit{Let $\phi:  {\mathbb{G}}_{m}^N \rightarrow  {\mathbb{G}}_{m}^N$ be a monomial map induced by $A$ with real eigenvalues, irreducible characteristic polynomial over $\mathbb{Q}$ and $\delta_{\phi}>1$, and let $P \in  {\mathbb{G}}_{m}^N(\bar{\mathbb{Q}})$ be a point with infinite orbit $\mathcal{O}_{\phi}(P)$. Suppose that the coordinates of $P$ are algebraic integers in $\mathbb{Q}(J(A),P)$. Then there is an effective computable positive constant $C$ depending on $J(A),  \mathbb{Q}(J(A),P)$ and $h(P)$ such that $ \hat{h}_{\phi_{A}}(P)>C$}.

%----------------------------------------------------------------------------------------
%	SECTION 2
%----------------------------------------------------------------------------------------

\section{Canonical heights and monomial maps}
Monomial maps are endomorphisms of the torus $\mathbb{G}_m^N$. They naturally induce self-rational maps of $\mathbb{P}^N$, by embedding $\mathbb{G}_m^N$ in $\mathbb{P}^N$. Before proving results, we state some definitions and results from sections 6 and 7 of [15].\newline

{\bf Definition 2.1:} \textit{Let $\phi : \mathbb{P}^N \dashrightarrow \mathbb{P}^N$ be a dominant rational map. The (first) dynamical degree of $\phi$ is the quantity} \begin{center} $\delta_{\phi}= \lim_{n \rightarrow \infty} \deg( \phi^n)^{\frac{1}{n}}$  \end{center}

{\bf Proposition 2.2:} \textit{The limit defining the dynamical degree exists and satisfies} \begin{center}
$\delta_{\phi}= \inf_{n \geq 1} \deg( \phi^n)^{\frac{1}{n}}$
\end{center}
\begin{proof}
See Proposition 7 of [15].
\end{proof}

Recall from the introduction that, for a point $P \in \mathbb{P}^N(\bar{\mathbb{Q}})$, the canonical height of $P$ with respect to $\phi$ is defined as follows: \newline \newline
{\bf Definition 2.3:} \textit{Let $\phi : \mathbb{P}^N \dashrightarrow \mathbb{P}^N$ be a dominant rational map with dynamical degree $\phi >1$. Let $P \in \mathbb{P}^N(\bar{\mathbb{Q}})_{\phi}$. The canonical height of $P$ is the quantity} 
\begin{center} 
$ \hat{h}_{\phi}(P)=\lim \sup_{n \rightarrow \infty} \dfrac{1}{n^{l_{\phi}} \delta_{\phi}^n} h(\phi^n(P))$ \end{center}

{\bf Example 2.4:} We assume $d>0$ and consider $\phi(P)=\phi(x_1,...,x_N)= P^{-d}:=(x_1^{-d},..., x_N^{-d})$. Therefore  $\phi^n(x_1,...,x_N)=(x_1^{-d^n},..., x_N^{-d^n})$ if $n$ is odd, and is $(x_1^{d^n},..., x_N^{d^n})$ if $n$ is even. Since the dynamical degree is $\delta_{\phi}=d$ and $l_{\phi}=0$, the sequence $\{ \dfrac{h(\phi^n(P)}{d^n} \}_{n=1}^{\infty}$ has two limit points, namely, $h(P)$ and $h(P^{-1})$. This implies that $ \hat{h}_{\phi}(P)= \max \{ h(P), h(P^{-1})\}.$ \newline \newline
Silverman proved some properties for this canonical height in [15], namely:

 (a) $0 \leq \hat{h}_{\phi}(P) \leq \infty.$

(b) $\hat{h}_{\phi}(\phi(P))=\delta_{\phi}\hat{h}_{\phi}(P).$

(c) If $P \in$ Preper$(\phi)$, then $\hat{h}_{\phi}(P)=0$
\newline \newline From now on, we focus our attention on monomial maps. A monomial map is an endomorphism of the torus $\mathbb{G}^N_m$. Embedding $\mathbb{G}^N_m$ in $\mathbb{P}^N$, monomial maps induce rational self-maps of $\mathbb{P}^N$ \newline \newline
{\bf Definition 2.5:}
\textit{We write ${\mbox{Mat}}_{\mathbb{Z}}^+$ for the set of $N$-by-$N$ matrices with integer coefficients and nonzero determinant. To each matrix $A= (a_{ij}) \in {\mbox{Mat}}_{\mathbb{Z}}^+$ we associate the monomial map $\phi_A: {\mathbb{G}}_{m}^N \rightarrow {\mathbb{G}}_{m}^N$ given by the formula}\newline

\textit{$\phi_A(X_1,...,X_N)=({X}_{1}^{a_{11}}{X}_{2}^{a_{12}}...{X}_{N}^{a_{1N}}, {X}_{1}^{a_{21}}{X}_{2}^{a_{22}}...{X}_{N}^{a_{2N}},..., {X}_{1}^{a_{N1}}{X}_{2}^{a_{N2}}...{X}_{N}^{a_{NN}})$.}\newline

 We call $\phi_A$ the monomial map associated to $A$, that induces a rational map $\phi_A:{\mathbb{P}}^{N} \dashrightarrow {\mathbb{P}}^{N}$. Again, we denote the spectral radius of $A$ by
\begin{center}
$\rho(A)= \max \{|\lambda|: \lambda \in \mathbb{C}  \mbox{ is an eigenvalue for} A\}. $
\end{center}
It is immediate from the definition that if $A, B \in {\mbox{Mat}}_{\mathbb{Z}}^+$ are matrices with associated monomial maps $\phi_A$ and $\phi_B$, then
\begin{center}
$\phi_{AB}(P)=(\phi_A \circ \phi_B)(P)$ and $\phi_{A+B}(P)= \phi_A(P). \phi_B(P).$
\end{center}
{\bf Definition 2.6:}
\textit{Let $A \in GL_N(\mathbb{Q})$. A Jordan subspace for $A$ is an $A$-invariant subspace of $\bar{\mathbb{Q}}^{N} $ corresponding to a single Jordan block of $A$. A Jordan subspace $V \subset\bar{\mathbb{Q}}^{N}$ with associated eigenvalue $\lambda$ is called a maximal Jordan subspace if $|\lambda|= \rho(A)$ and if the dimension of $V$ is maximal among the Jordan subspaces whose eigenvalue has magnitude equal to $\rho(A)$. We set}\newline \newline
\textit{$r(A)=$ number of maximal Jordan subspaces, ~~~~~~~~~~~~~~~~~~~~~~~~~~~~~~~~~~~~~~~~~~~~~~~~~~~\newline \newline
$\bar{r}(A)= \# \{ \sigma(V) : V \mbox{is a maximal Jordan subspace for} ~ A ~ \mbox{and} ~ \sigma \in $ Gal$(\bar{\mathbb{Q}}/ \mathbb{Q}) \},$} \newline \newline
$l(A)= \dim($\textit{any maximal Jordan subspace}$)-1$. ~~~~~~~~~~~~~~~~~~~~~~~~~~~~~~~~~~~~~~~~~~~~~~~~\newline

Thus, $\bar{r}(A)$ is the number of distinct $\bar{\mathbb{Q}}$-subspaces of $\bar{\mathbb{Q}}^{N}$ that are Galois conjugate to a maximal Jordan subspace of $A$, and so $\bar{r}(A) \geq r(A) \geq 1$, since $A$ always has at least one maximal Jordan subspace.\newline

{\bf Definition 2.7:}
\textit{ Let $G$ be an algebraic subgroup of ${\mathbb{G}}_{m}^N$. We write $G(\bar{\mathbb{Q}})^{\mbox{div}}$ for the divisible hull of $G(\bar{\mathbb{Q}})$,}
\begin{center}
\textit{$G(\bar{\mathbb{Q}})^{\mbox{div}}= \{(\alpha_1,..., \alpha_N) \in {\mathbb{G}}_{m}^N(\bar{\mathbb{Q}})  :   ({\alpha}_{1}^{n},..., {\alpha}_{N}^{n}) \in G(\bar{\mathbb{Q}})$ for some $n \geq 1 \}$.}
\end{center}
\textit{Equivalently, $G(\bar{\mathbb{Q}})^{\mbox{div}}$ is the set of translates of $G(\bar{\mathbb{Q}})$ by points in ${\mathbb{G}}_{m}^N(\bar{\mathbb{Q}})_{\mbox{tors}}$. }\newline
\newline 
 {\bf Theorem 2.8}: \textit{Let $A \in {\mbox{Mat}}_{\mathbb{Z}}^+$ be a matrix with associated monomial map $\phi_A$.\newline
(a) $\rho (A) \geq 1$.
\newline(b) $\rho(A)=1$ if and only if all the eigenvalues of $A$ are roots of unity.
\newline(c) (Hasselblatt, Prop) $ \delta_{\phi_A}= \rho(A)$.} \newline
(d) [Lin, 11, Theorem 6.2]  \textit{$l_{\phi_A}=l(A).$} 

 \begin{proof}
See Proposition 21 and Theorem 24 of [15].
\end{proof}

We remark why points with canonical height equal to zero lie on a non-dense subgroup of the torus, remaking here the proof done by Silverman. \newline

{\bf Theorem 2.9:} [Silverman, 15] \textit{Let $A  \in {\mbox{Mat}}_{\mathbb{Z}}^+$ be a matrix whose associated monomial map $\phi_A$ has dynamical degree $\delta_{\phi_A}>1$. Then there is an algebraic subgroup $G \subset {\mathbb{G}}_{m}^N$ of dimension}
\begin{center}
$\dim G \geq N - \bar{r}(B)$
\end{center} \textit{such that}
\begin{center}
$\{ P \in {\mathbb{G}}_{m}^N(\bar{\mathbb{Q}}) ; \hat{h}_{\phi_A}(P)=0 \} \subset G(\bar{\mathbb{Q}})^{\mbox{div}}$.
\end{center}
\begin{proof} We reproduce Silverman's proof here. Making $\phi= \phi_A ,\rho= \rho(A)= \delta_{\phi_A}, l=l(A)=l_{\phi_A},$ $K$ a number field which is Galois over $\mathbb{Q}$ and contains the eigenvalues of $A$ and the coordinates of $P$, $M_K$ its set of places, and $||.||_v$ the absolute value associated to the place $v \in M_K$. We consider $P \in \mathbb{G}_m^N(\bar{\mathbb{Q}})$, such that $\hat{h}_{\phi}(P)=0$. The sequence of matrices $\{\dfrac{A^n }{n^l \rho^n} : n \geq 0 \}$ has an accumulation point in $\mathbb{C}^{N^2}$. We denote \begin{center} $\max^+((u_1,...,u_N))= \max(0,u_1,...u_N)$ for all $(u_1,...,u_N) \in \mathbb{R}^{N}$,\end{center}  and call \begin{center}$\log||(y_1,...,y_N)||_v:= (\log||y_1||_v,...,\log||y_N||_v)$ for all $(y_1,...,y_N) \in \overline{\mathbb{Q}}^N$. \end{center} Then $\log ||\phi_A^n(P)||_v = A^n \log || P ||_v.$ and choose $\mathcal{N} \subset \mathbb{N}$ such that

\begin{center} $\hat{h}_{\phi}(P)= \lim \sup_{n \rightarrow \infty} \dfrac{h(\phi^n(P))}{n^l \rho^n}=\lim_{n \in \mathcal{N}} \dfrac{h(\phi^n(P))}{n^l \rho^n}= 0$ 
\end{center} and \begin{center} $B := \lim_{n \in \mathcal{N}} \dfrac{A^n }{n^l \rho^n}. \:\:\:\:\:(1)$ \end{center} Thus  \begin{center}
$0= \sum_{v \in M_K} \max^{+}(\lim_{n \in \mathcal{N}} \dfrac{A^n }{n^l \rho^n} \log || P ||_v)= \sum_{v \in M_K} \max^{+}( B \log || P ||_v).$
\end{center} Then  $\log || P ||_v \in \ker_{\mathbb{C}}(B)$ for all $v \in M_K.$ \:\:\:(2)

Let also $V \subset \bar{\mathbb{Q}}^N$ be a Jordan subspace for $A$, which means that $A_{|V}$ is a matrix in Jordan form. If $V$ is not a maximal Jordan subspace, one can check that $\lim_{n \in \mathcal{N}} {(\dfrac{A^n}{n^l \rho^n})}{|V}=B_{|V}=0$.

If $V$ is maximal and $\mathcal{V}= \{v_1,...,v_t \}$ a basis of $V$ used to put $A_{|V}$ into Jordan normal form. Then $Bv_1=...=Bv_{t-1}=0$ and $B v_t=\dfrac{\xi}{\rho^l l!} v_1$ for some $\xi \in \mathbb{C}$ with $|\xi|=1.$ \:\:\:(3)

We can write $\bar{\mathbb{Q}}^N=V_1+...+V_r + Z$, where $V_1,...,V_r$ are the distinct maximal Jordan subspaces and $Z$ is the sum of all the other Jordan subspaces for $A$. Moreover, for each $i$ we let $W_i \subset V_i$ be the $A$-invariant codimension 1 $\bar{\mathbb{Q}}$-subspace of $V_i$ satisfying  $\ker(B|_{V_i \otimes_{\bar{\mathbb{Q}}} \mathbb{C}})=W_i \otimes_{\bar{\mathbb{Q}}} \mathbb{C}.$ Hence \begin{center} 
$\ker_{\mathbb{C}}(B)= U \otimes_{\bar{\mathbb{Q}}} \mathbb{C},$ where $U= W_1+...+W_r + Z \subset \bar{\mathbb{Q}}^N, \dim_{\bar{\mathbb{Q}}} U= N-r.\:\:\:(4)$ 
\end{center}
For $W \subset K^N$, and $F \subset$ a subfield, denote the subspace of $F^N$ that is orthogonal to $W$ by 
\begin{center}
Perp$_{F}(W)= \{b \in F^N: b.w=0$ for all $w \in W\}$.
\end{center}
Hence \newline \newline 
Perp$_{\mathbb{C}}(\log ||P||_v) \supset$ Perp$_{\mathbb{C}}(\ker_{\mathbb{C}}(B))=$Perp$_{\mathbb{C}}(U \otimes_{\bar{\mathbb{Q}}} \mathbb{C})=$Perp$_{\bar{\mathbb{Q}}}(U) \otimes_{\bar{\mathbb{Q}}} \mathbb{C} \supset$ Perp$_{\bar{\mathbb{Q}}}(U).$  \newline \newline

Then there is unique $\mathbb{Q}$-vector space $Y \subset \mathbb{Q}$ depending only on $A$ such that \begin{center}   Perp$_{\bar{\mathbb{Q}}}(\log ||P||_v) \supset Y$.          \end{center}

Defining the lattice $L:= Y \cap \mathbb{Z}^N$, and the assocated algebraic subgroup \begin{center}$G_L= \bigcap_{(e_1,...,e_N) \in L} \{ X_1^{e_1}...X_N^{e_N}=1 \} \subset \mathbb{G}^N_m,$
\end{center} we have that \begin{center}
$\{ P \in {\mathbb{G}}_{m}^N(\bar{\mathbb{Q}}) ; \hat{h}_{\phi_A}(P)=0 \} \subset G_L(\bar{\mathbb{Q}})^{\mbox{div}}$. 
\end{center}
The missing part is a calculation on the dimension of $Y$ and its definition.

\end{proof}

{\bf Corollary 2.10:} [Silverman, 15] \textit{Let $\phi_A \subset \mbox{End}({\mathbb{G}}_{m}^N)$ as above with $\delta_{\phi_A}>1$, and let $P$ be a point with $ \hat{h}_{\phi_A}(P)=0$. Then there is a proper algebraic subgroup $G$ of ${\mathbb{G}}_{m}^N$ with $\mathcal{O}_{\phi_A}(P) \subset G$. In particular, the orbit $\mathcal{O}_{\phi_A}(P)$ is not Zariski dense in ${\mathbb{G}}_{m}^N$.}
\begin{proof} As $ \hat{h}_{\phi_A}(P)=0$, it follows by height properties that $\mathcal{O}_{\phi_A}(P) $ is a subset of $ \{ P \in {\mathbb{G}}_{m}^N(\bar{\mathbb{Q}}) ; \hat{h}_{\phi_A}(P)=0 \}$, which is contained in $G_L(\bar{\mathbb{Q}})^{\mbox{div}}$ by the previous Theorem. Let $d$ be the number of roots of unity in $K$. Then $\mathcal{O}_{\phi_A}(P)^d \subset G_L(K)$, and therefore $\mathcal{O}_{\phi_A}(P) \subset G_{dL}(K)$.
\end{proof}

{\bf Corollary 2.11:} [Silverman, 15] \textit{Let $\phi_A \subset \mbox{End}({\mathbb{G}}_{m}^N)$ be a monomial map induced by $A \in {\mbox{Mat}}_{\mathbb{Z}}^+$ with $\delta_{\phi_A}>1$, whose characteristic polynomial is irreducible over $\mathbb{Q}$. Let $P \in {\mathbb{G}}_{m}^N(\bar{\mathbb{Q}})$. Then}
\begin{center}
$ \hat{h}_{\phi_A}(P)=0 \iff \# \mathcal{O}_{\phi_A}(P) < + \infty$. 
\end{center}
\begin{proof}  It is already known that $\hat{h}_{\phi_A}(P)=0$ if $P \in$ Preper$(\phi)$. Conversely, since the characteristic polynomial is irreducible, the eigenvalues $\lambda_1,..., \lambda_N$ of $A$ are all distinct, the Jordan subspaces are all $1$-dimensional and $\bar{r}(A)=N$. Hence the algebraic group $G$ from Theorem 2.9 has dimension $0$, and $G(\bar{\mathbb{Q}})^{\mbox{div}} = {\mathbb{G}}_{m}^N(\bar{\mathbb{Q}})_{\mbox{tors}}$. Thus $\hat{h}_{\phi_A}(P)=0$ implies that $P 
\in {\mathbb{G}}_{m}^N(\bar{\mathbb{Q}})_{\mbox{tors}}$, and therefore $P$ is preperiodic.
\end{proof}

\section{Effective bounds for the canonical height of non-periodic points}

 Silverman observes in [15, Remark 30] that it should be possible to use an effective form of Baker's theorem to prove effective versions of the results above, with an effective computable constant $C=C(A,h(P))>0$ such that 
\begin{center}
$\mathcal{O}_{\phi_{A}}(P)$ Zariski dense $\implies \hat{h}_{\phi_{A}}(P)>C. $
\end{center}
In this section we work out  the cases for that $A$ has real Jordan form, i.e, all of its eigenvalues are real, obtaining a constant $C$ that depends also on the number field extension where $P$ and the Jordan form of $A$ are defined, which we can call $\mathbb{Q}(J(A),P)$. For points whose coordinates are not necessarily integral over this field, although integral after being multiplyied by an integer $\zeta$, the desired constant will also depend on this $\zeta$ for instance, which is related with the arithmetic complexity of the point. 

 In order to prove our results, we will make use of an improvement of an effective classical Baker's Theorem that is due to P. Philippon and M. Waldschmidt, and can be seen for example in [3, chapter 18, Theorem 1.1]. \newline

{\bf Theorem 3.1:} [P. Philippon and M. Waldschmidt]
\textit{Let $\alpha_1,..., \alpha_n$ be non-zero algebraic numbers that are different from 1, and let $\beta_1,..., \beta_n$  algebraic numbers not all zero, and $\log \alpha_1,..., \log \alpha_n$ logarithmic representatives such that $\pi i, \log \alpha_1,..., \log \alpha_n$ are linearly independent over $\bar{\mathbb{Q}}$ and $\Lambda:= \beta_1 \log \alpha_1 + \beta_2 \log \alpha_2 +...+\beta_n \log \alpha_n$ does not vanish. }

\textit{Let $D$ be a positive integer, $A, A_1, A_2,..., A_n$ be positive real numbers, and $B$ satisfying}
\begin{center}
$D \geq [\mathbb{Q}(\alpha_1,...,\alpha_n, \beta_1,..., \beta_n) : \mathbb{Q}]$
\end{center}
\begin{center}
$A_j \geq \max \{ H(\alpha_j), \exp|\log \alpha_j|, e^n \}, ~~~~~ 1 \leq j \leq n$
\end{center}
\begin{center}
$A:= \max \{A_1,..., A_n, e^e \}$
\end{center}
\begin{center}
$B:= \max \{H(\beta_j) ; 1 \leq j \leq n \}$
\end{center}
\textit{Then}

\begin{center}
$| \Lambda| \geq e^{-U},$
\end{center}

\textit{where} \begin{center}$U=-C_{11}(n).D^{n+2}.\log A_1 ... \log A_n. (\log B + \log \log A)$, \end{center} \begin{center}  $C_{11}(n) \geq 2^{8n +53}.n^{2n}$, \end{center} \textit{ and $H(\alpha)$ denotes the maximun of absolute values of the coefficients of the minimum polynomial of $\alpha$ over $\mathbb{Q}$.}

\begin{proof}  Can be found in [3, chapter 18, section 4].

\end{proof}

We remind that, according to the proof of  Theorem 6.2 of [11], it is true for monomial maps $\phi$ that $\inf_{n \geq 0} \dfrac{\deg (\phi^n)}{n^{l_{\phi}}\delta_{\phi}^n} >0.$

For a matrix $A$, and a point $P$ on a projective space, we denote by $J(A)$ a matrix that changes the canonical basis for another one that makes $A$ to assume its Jordan normal form, and by $\mathbb{Q}(J(A),P)$ a number field containing the coefficients of $P$ and $J(A)$ that is Galois over $\mathbb{Q}$. \newline \newline
{\bf Theorem 3.2:}
\textit{Let $\phi:  {\mathbb{G}}_{m}^N \rightarrow  {\mathbb{G}}_{m}^N$ be a monomial map induced by $A$ with real eigenvalues, $l(A) \geq 1$ and $\delta_{\phi}>1$, and let $P \in  {\mathbb{G}}_{m}^N(\bar{\mathbb{Q}})$ be a point with orbit $\mathcal{O}_{\phi}(P)$ Zariski dense. Suppose $\zeta \in \mathbb{Z}$ satisfying that $\zeta P$ has algebraic integer coordinates, and $C_0=\inf_n \dfrac{\deg(\phi^n)}{n^l \rho^n} >0$ are both well-known. Then there are an effective computable positive constant C depending on $\mathbb{Q}(J(A),P), A$ (more precisely on the heights of the coordinates of $J(A)$), and $h(P)$, and  an effective function $m(\zeta , C)$ such that $ \hat{h}_{\phi_{A}}(P)>C+ C_0 m(\zeta, C )$}.

\begin{proof} We denote \begin{center} $\max^+((u_1,...,u_N))= \max(0,u_1,...u_N)$ for all $(u_1,...,u_N) \in \mathbb{R}^{N}$,\end{center}  and call \begin{center}$\log||(y_1,...,y_N)||_v:= (\log||y_1||_v,...,\log||y_N||_v)$ for all $(y_1,...,y_N) \in \overline{\mathbb{Q}}^N$, \end{center} and we call $P=(x_1,...,x_N)$.

Starting from $A \in {\mbox{Mat}}_N^+(\mathbb{Z})$, we can write $\bar{\mathbb{Q}}^N=V_1+...+V_r + Z$, where $V_1,...,V_r$ are the distinct maximal Jordan subspaces and $Z:= V_{r+1}+...+V_t$ is the sum of all the other Jordan subspaces for $A$, as in the proof of Theorem 2.9.

 In [14], we have an effective algorithm for finding a basis where $A$ has Jordan normal form defined over some algebraic extension of $\mathbb{Q}$. In other words, for each $1 \leq i \leq t$, we can find the basis $\{ {v}_1^{(i)},..., {v}_{t_i}^{(i)}\}$ used to put $A_{|V_i}$ in the Jordan normal form, with the canonical coordinates denoted by ${v}_j^{(i)}=({a}_{1j}^{(i)},..., {a}_{Nj}^{(i)})$, effectively computable, defined over some algebraic extension of $\mathbb{Q}$, that we can suppose to be Galois, and that we call $K$. 

Denoting again $l:=l_{\phi_A}=l(A)$, $r(A)=r$, and $\rho:= \rho(A)$, the spectral radius of $A$, we have by (1) of Theorem 2.9 that there exists an infinite subset $\mathcal{N} \subset \mathbb{N}$ and $B \in {\mbox{Mat}}_N(\mathbb{R})$ such that the limit $B=\lim_{n \in \mathcal{N}} \dfrac{A^n}{{n}^l {\rho}^n}$ is satisfied.

To suppose $\mathcal{O}_{\phi_{A}}(P)$ Zariski dense implies that $ \hat{h}_{\phi_{A}}(P)>0$ by Theorem 2.9. Extending $K$, we can suppose $P=(x_1,...,x_N)$ defined over $K$. By (2) on the proof of Theorem 2.9  we have that \begin{center}$\hat{h}_{\phi_{A}}(P)=\sum_{v \in M_K}\max^+(B\log||P||_v) > 0$ \end{center} and \begin{center}$ \log||P||_v \notin \ker_{\mathbb{C}}(B) ~ \forall ~ v$ with nonzero $v$-component in the above sum. \end{center} By (4) of Theorem 2.9 we can see that $ \log||P||_v \notin \ker_{\mathbb{C}}(B)$ implies that \begin{center} $ \log||P||_v \in (\mathbb{C}{v}_{t_1}^{(1)} \bigoplus ... \bigoplus \mathbb{C}{v}_{t_r}^{(r)} \bigoplus \ker_{\mathbb{C}}(B)) -  \ker_{\mathbb{C}}(B)$. \end{center} For any $v \in M_K$, we can make an effective change of basis using Cramer rule to obtain $\log||P||_v$ in the basis where $A$ is in Jordan form, in other words, to obtain $c_{it_i,v}(P) \in \mathbb{C}, i \leq r, b \in  \ker_{\mathbb{C}}(B)$, such that $\log||P||_v = \sum_{i \leq r}c_{it_i,v}(P){v}_{t_i}^{(i)} + b$, and so \begin{center} $\log||P||_v=(\sum_{i\leq r}c_{it_i,v}(P){a}_{1t_i}^{(i)},...,\sum_{i\leq r}c_{it_i,v}(P) {a}_{Nt_i}^{(i)}) +b$.
\end{center} 
Effectively, let $J(A)$ be the $(N\times N)$-matrix $({a}_{il}^{(j)})$ with lines indexed by $1 \leq i \leq N$ and columns indexed  by $(j,l); 1 \leq j \leq t, 1 \leq l \leq t_j$ in lexicographic order. So $c_{it_i,v}(P)$ are coodinates of the vector solution $z$ for the linear system $J(A).z=\log||P||_v$. Using Cramer Rule we see that these solutions have the form $c_{it_i,v}(P)=\sum _{j\leq N}d_{ij,v}(P) \log||x_j||_v$, for $d_{ij,v}(P) \in K$ effectively computable depending only on $A$.

 From (3) of Theorem 2.9, we see that $B{v}_{t_i}^{(i)} =\frac{\xi_i}{{\rho}^l l!}{v}_{1}^{(i)}$ for $\xi_i \in \mathbb{C}, |\xi_i|=1$ for $1 \leq i \leq r$, and $ \xi_i = \lim_n \dfrac{\lambda^n}{\rho^n}$ for some eigenvalue $\lambda$ of $A$. The eigenvalues of $A$ are real by hypothesis, therefore $\xi_i$ is equal to $1$ or $-1$.  Then for this case we have that\newline \newline
 $B\log ||P||_v = \sum_{i \leq r}c_{it_i,v}(P)B{v}_{t_i}^{(i)}=\sum_{i \leq r}c_{it_i,v}(P)\dfrac{\xi_i}{{\rho}^l l!}{v}_{1}^{(i)} ~~~~~~~~~~~~~~~~~~~~~~~~~~~~~~~~\newline 
=\dfrac{1}{{\rho}^l l!}(\sum_{i \leq r}c_{it_i,v}(P) {a}_{11}^{(i)}  \xi_i,...,\sum_{i \leq r}c_{it_i,v}(P) {a}_{N1}^{(i)}  \xi_i )=~~~~~~~~~~~~~~~~~~~~~~~~~~~~~~~~~~~~~~~\newline
 \dfrac{1}{{\rho}^l l!}(\sum_j(\sum_{i \leq r}{a}_{11}^{(i)}  \xi_id_{ij,v}(P))\log||x_j||_v,..., \sum_j(\sum_{i \leq r}{a}_{N1}^{(i)}  \xi_id_{ij,v}(P))\log||x_j||_v)$,\newline \newline which yields \newline \newline
$   \hat{h}_{\phi_{A}}(P)=\sum_{v \in M_K}\max^+(B\log||P||_v) ~~~~~~~~~~~~~~~~~~~~~~~~~~~~~~~~~~~~~~~~~~~~~~~\newline
=\dfrac{1}{{\rho}^l l!}\sum_{v \in M_K}\max \{0,\sum_j(\sum_{i \leq r}{a}_{t1}^{(i)}  \xi_id_{ij,v}(P))\log||x_j||_v  ;t \leq N \} 
> 0$.\newline

Multiplying the coordinates of $P$ by an adequate algebraic integer $\zeta \in K$, one can suppose for instance that each coordinate of $P$ is an algebraic integer of $K$. The canonical height of the point is independent of this change, and following analogous calculations yields that:  \newline \newline
$   \hat{h}_{\phi_{A}}(P) =\newline\sum_{v \in M_K}\max \{ \lim_{n \in \mathcal{N}} \dfrac{\deg(\phi^n)}{n^l \rho^n}\log ||\zeta||_v,\dfrac{1}{{\rho}^l l!} \sum_j(\sum_{i \leq r}{a}_{t1}^{(i)}  \xi_id_{ij,v}(P))\log||x_j||_v  ;t \leq N \} \newline
> 0$, \newline \newline  and then \newline \newline $ \hat{h}_{\phi_{A}}(P) \newline \geq  \sum_{v \in M_K} \max \{ C_0 \log ||\zeta||_v,\dfrac{1}{{\rho}^l l!} \sum_j(\sum_{i \leq r}{a}_{t1}^{(i)}  \xi_id_{ij,v}(P))\log||x_j||_v  ;t \leq N \} >0$
\newline \newline for $C_0= \inf_n \dfrac{\deg(\phi^n)}{n^l \rho^n} >0$. Now consider $S\in M_K$ the finite set of places $v$ such that \begin{center}$\max \{\sum_j(\sum_{i \leq r}{a}_{t1}^{(i)}  \xi_id_{ij,v}(P))\log||x_j||_v  ;t \leq N \}> C_0 \log ||\zeta||_v$, \end{center} and $R = M_K - S$  the set of all the other places. For example $S \subset T := \{v ; ||x_j||_v \neq 1 ~\text{ for some}~ j\}$ finite. Denote by $S^0, S^{\infty}, T^0, T^{\infty}, M_K^0, M_K^{\infty}$ the sets of non-archimedean places and the sets of archimedean places for each of the sets $S, T, M_K$.

Suppose that for $v \in S$, the maximum amongst the coordinates is achieved for $t=t_v$, and denote $D_{ij,v}:={a}_{t_v1}^{(i)}.d_{ij,v}(P)$. We use the notation $N(v)$ for the norm of the ideal in $K$ corresponding to the place $v$. 

We aim to make use of the fact that $\{ \log N(v) ; v \in S_0 \}$ is linearly independent over $\overline{\mathbb{Q}}$ to apply Theorem 3.1.
 We start computing\newline \newline
$   \hat{h}_{\phi_{A}}(P)=\dfrac{1}{{\rho}^l l!}\sum_{v \in S} \sum_j(\sum_{i \leq r}{a}_{t_v1}^{(i)}  \xi_id_{ij,v}(P))\log||x_j||_v + C_0\sum_{v \in R} \log||\zeta||_v~~~~~~~~ ~~~~~~~~~~~\newline
=\dfrac{1}{{\rho}^l l!}\sum_{v \in S} \sum_j(\sum_{i \leq r} D_{ij,v} \xi_i)\log||x_j||_v + C_0\sum_{v \in R} \log||\zeta||_v ~~~~~~~~~~~~~~~~~~~~~~~~~~~~~~~~~~~~~~~~~~~~\newline
 =\dfrac{1}{{\rho}^l l![K: \mathbb{Q}]}\sum_{v \in S^0} [\sum_{i,j} D_{ij,v} \xi_i(-v(x_j))+\sum_{i,j,u \in S^{\infty}} D_{ij,u} \xi_i.v_{u,j}].\log N(v) + C_0\sum_{v \in R} \log||\zeta||_v$,\newline \newline where $v_{u,j} \in \mathbb{Q}$ and $\sum_{u \in M_K^{\infty}}v_{u,j}=v(x_j)$ for $v \in S^0$ a non-archimedean place. $K$ is a Galois extension of $\mathbb{Q}$, thus $v_{u,j} = v_{w,j}=$ for all $u, w \in S^{\infty}$.

If the sum $\dfrac{1}{{\rho}^l l!}\sum_{v \in S} \sum_j(\sum_{i \leq r}{a}_{t_v1}^{(i)}  \xi_id_{ij,v}(P))\log||x_j||_v $ is equal to $0$, one could make$C=0$ for instance. Otherwise, it follows that 
\begin{center}
$ \sum_{v \in S^0} [\sum_{i,j} D_{ij,v} \xi_i(-v(x_j))+\sum_{i,j,u \in S^{\infty}} D_{ij,u} \xi_i.v_{u,j}].\log N(v)  \neq 0$
\end{center}
We denote $H_{\bar{\mathbb{Q}}}:= \exp h$. Using the facts that $H_{\bar{\mathbb{Q}}}(a_1 +...+ a_n) \leq nH_{\bar{\mathbb{Q}}}(a_1)...H_{\bar{\mathbb{Q}}}(a_n)$, and that $H(\alpha) \leq (2.H_{\bar{\mathbb{Q}}}(\alpha))^{[K: \mathbb{Q}]}$[19, first inequation of page 77 and Lemma 3.11], \newline  where $H_{\bar{\mathbb{Q}}}:=\exp{h}$ is the Weil multiplicative height, we have that \newline \newline
$H(\sum_{i,j} D_{ij,v} \xi_i(-v(x_j))+\sum_{i,j,u \in S^{\infty}} D_{ij,u} \xi_i.v_{u,j}) ~~~~~~~~~~~~~~~~~~~~~~~~~~~~~~~~~~~~~~~~~~~~~~\newline \newline
 \leq [2. H_{\bar{\mathbb{Q}}}(\sum_{i,j} D_{ij,v} \xi_i(-v(x_j))+\sum_{i,j,u \in S^{\infty}} D_{ij,u} \xi_i.v_{u,j})]^{[K:\mathbb{Q}]} ~~~~~~~~~~~~~~~~~~~~~~~~~~~~\newline \newline 
 \leq [2Nr([K:\mathbb{Q}]+1)  \prod_{i,j,u \in S^{\infty}}H_{\bar{\mathbb{Q}}}( D_{ij,v} \xi_i(-v(x_j)) H_{\bar{\mathbb{Q}}}(D_{ij,u} \xi_i.v_{u,j})]^{[K:\mathbb{Q}]} 
~~~~~~~~~~~~~~~~~\newline \newline
\leq (4Nr.[K:\mathbb{Q}].\max \{[K : \mathbb{Q}],\max_{j}|v(x_j)|\}^{2Nr{[K:\mathbb{Q}]}}.\max_{i,j,w \in T}H_{\bar{\mathbb{Q}}}(D_{ij,w})^{2Nr[K: \mathbb{Q}]})^{[K:\mathbb{Q}]} ~ \newline \newline \leq
 \{4Nr.[K:\mathbb{Q}]^2.\max_{j}|v(x_j)|.(N-1)!\max_{i,j,l}H_{\bar{\mathbb{Q}}}(a_{il}^{(j)})H_{\bar{\mathbb{Q}}}(\det J(A))\}^{2N^2r[K:\mathbb{Q}]^2} \newline \newline \leq
 \{4Nr.[K:\mathbb{Q}]^2.\max_{j}|v(x_j)|.(N-1)! N! [\max_{i,j,l}H_{\bar{\mathbb{Q}}}(a_{il}^{(j)})]^2\}^{2N^2r[K:\mathbb{Q}]^2}.$\newline \newline 
 Then we have by Theorem 3.1 that 
\newline \newline
$\hat{h}_{\phi_{A}}(P) \geq \dfrac{1}{{\rho}^l l! [K: \mathbb{Q}]} \exp(-{E}.\prod_{v \in T^0}\log A_v.[\log D +  \max_{v \in T^0}\log \log A_v]) + C_0\sum_{v \in R} \log||\zeta||_v$,~~~~~~~~~~
\newline \newline where 

$A_v::=\max \{ \exp |\log N(v)|, e^{\# T^0}, e^e \} \leq N(v)^{6.\#T^0}$ for every $v \in T^0$, \newline

${E}:=C_{11}(\#T^0).[K:\mathbb{Q}]^{\#T^0 +2}, $\newline

 $D:=\max_{v \in T^{0}}  \{4Nr.[K:\mathbb{Q}]^2.\max_{j}|v(x_j)|.(N-1)! N! [\max_{i,j,l}H_{\bar{\mathbb{Q}}}(a_{il}^{(j)})]^2\}^{2N^2r[K:\mathbb{Q}]^2}$,\newline

  for $T^0= \{v \in M_K ; ||x_j||_v \neq 1 ~\text{ for some}~ j\}$, $C_{11}(n) = 2^{8n +53}.n^{2n}$. \newline

We denote $h_K(P)= [K:\mathbb{Q}]. h(P)$, and use the product formula and the fact that the $x_j'$s are algebraic integers  to conclude that \newline \newline 
$\sum_{v \in T^0} \log \max \{1, N(v)^{|v(x_j)|} \} =  \sum_{v \in T^0} \log N(v)^{v(x_j)}=\sum_{v \in T^{\infty}} \log ||x_j||_v  \leq \sum_{v \in T^{\infty}} \log \max \{1, ||x_j||_v \} \newline \newline \leq h_K(x_j) \leq h_K(P),$ and then  \newline \newline
$\#T^0=\sum_{v \in T^0} 1\leq  \sum_j \sum_{v \in T^0} \log \max \{1, N(v)^{|v(x_j)|} \} \leq 1 + Nh_K(P), \newline \newline
|v(x_j)|\leq  \prod_j \prod_{v \in T^0} \max \{1, N(v)^{|v(x_j)|} \} \leq (H_K(P))^N,$ where $H_K(P)= \exp(h_K(P)) \newline \newline
\log N(v) \leq \log N(v)^{|v(x_j)|} \leq \sum_j \sum_{v \in T^0} \log \max \{1, N(v)^{|v(x_j)|} \} \leq  Nh_K(P),$ \newline \newline

 Making $C$ equal to \begin{center} $\dfrac{1}{2{\rho}^l l![K: \mathbb{Q}]} \exp(-{E}^{\prime}.(\log A^{\prime})^{(1+Nh_K(P))}.[\log D^{\prime} +\log \log A^{\prime}])$, \end{center}
where \newline

$A^{\prime}::= e^{( Nh_K(P))(6.(1 + Nh_K(P)))}$,\newline

${E}^{\prime}:=C_{11}(1 + Nh_K(P)).[K:\mathbb{Q}]^{(3 + Nh_K(P))}, $\newline

 $D^{\prime}:= \{4Nr.[K:\mathbb{Q}]^2.(H_K(P))^N.(N-1)!N![\max_{i,j,l}H_{\bar{\mathbb{Q}}}(a_{il}^{(j)})]^2 \}^{2N^2r[K:\mathbb{Q}]^2}$,\newline \newline
 we have that \begin{center}$  \hat{h}_{\phi_{A}}(P)>C  + C_0\sum_{v \in R} \log||\zeta||_v $\end{center} The function $m$ to be defined will arise then from the factorization of $|\zeta| \in \mathbb{N} - 0$. Namely, let us suppose that the factorization of $|\zeta|$ with prime numbers is given by $|\zeta|= p_1^{e_1}...p_n^{e_n}$, which can be found algorithmically, then we define:
\begin{center}
$m(\zeta,C):= \min \{m= \sum_{1 \leq i \leq n} \log(p_1^{t_1}...p_n^{t_n}) | t_i= -e_i, 0, $ or $ e_i ,$ and $ m+C>0 \}.$
\end{center} Thus it follows that \begin{center}$  \hat{h}_{\phi_{A}}(P)>C  + C_0 m( \zeta , C)$\end{center} as wanted, noting that $h_K(P)$ is independent of the projective representant for $P$, so there was no problem in multiplying $P$ by $\zeta$.
\end{proof} 

{\bf Corollary 3.3:}
\textit{Let $\phi:  {\mathbb{G}}_{m}^N \rightarrow  {\mathbb{G}}_{m}^N$ be a monomial map induced by $A$ with real eigenvalues, irreducible characteristic polynomial over $\mathbb{Q}$ and $\delta_{\phi}>1$, and let $P \in  {\mathbb{G}}_{m}^N(\bar{\mathbb{Q}})$ be a point with infinite orbit $\mathcal{O}_{\phi}(P)$. Suppose $\zeta \in \mathbb{Z}$ satisfying that $\zeta P$ has algebraic integer coordinates, and $C_0=\inf_n \dfrac{\deg(\phi^n)}{n^l \rho^n} >0$ are both well-known. Then there is an effective computable positive constant C depending on $\mathbb{Q}(J(A),P), A$ (more precisely on the heights of the coordinates of $J(A)$), and $h(P)$,   and  an effective function $m(\zeta, C )$ such that $ \hat{h}_{\phi_{A}}(P)>C+ C_0 m(\zeta, C )$}.
\begin{proof}
By Corollary 2.11, we must have $\hat{h}_{\phi_A}(Q)>0$. Now we follow the proof of theorem 4.2, with some reductions and changes. Again we have $\bar{\mathbb{Q}}^N=V_1+...+V_r + Z$, where $V_1,...,V_r$ are the distinct maximal Jordan subspaces and $Z:= V_{r+1}+...+V_t$ is the sum of all the other Jordan subspaces for $A$. We assume $v_i:=(a_1^{(i)},..., a_N^{(i)})$ puts $A_{|V_i}$ in the Jordan normal form. Again $B:=\lim_{n \in \mathcal{N}} \dfrac{A^n}{{n}^l {\rho}^n}=\lim_{n \in \mathcal{N}} \dfrac{A^n}{{\rho}^n}$, but now (3) in the proof of Theorem 2.9 shows that \begin{center} $Bv_i= \lim_n \dfrac{\lambda^n}{\rho^n}v_i= \xi_iv_i, \xi_i \in \{-1, 1\} ~ \forall ~ 1 \leq i \leq r$. \end{center} Similarly to Theorem 3.2, we have that $\log||P||_v = \sum_{i \leq r}c_{i,v}(P){v}_{i} + b$  for $c_{i,v}(P) \in \mathbb{C}, i \leq r, b \in  \ker_{\mathbb{C}}(B)$, effectively computable and $K$ a number field such that the Jordan normal form of $A$ and the point $P$ are defined over it.   So \begin{center} $\log||P||_v=(\sum_{i\leq r}c_{i,v}(P){a}_{1}^{(i)},...,\sum_{i\leq r}c_{i,v}(P) {a}_{N}^{(i)}) +b$. \end{center}
Again $c_{i,v}(P)=\sum _{j\leq N}d_{ij,v}(P) \log||x_j||_v$, for $d_{ij,v}(P) \in K$ effectively computable by Cramer rule, depending only on $A$ implies that \newline \newline
 $B\log ||P||_v=  (\sum_j(\sum_{i \leq r}{a}_{1}^{(i)}  \xi_id_{ij,v}(P))\log||x_j||_v,..., \sum_j(\sum_{i \leq r}{a}_{N}^{(i)}  \xi_id_{ij,v}(P))\log||x_j||_v)$. \newline \newline
For some $S \subset T := \{v ; ||x_j||_v \neq 1 ~\text{ for some}~ j\}$ finite, $t_v \in \{1,...,r\}$ , denoting $D_{ij,v}:={a}_{t_v}^{(i)}.d_{ij,v}(P)$, and \newline

$A::=e^{( Nh_K(P))(6.(1+ Nh_K(P)))}$,\newline

${E}:=C_{11}(1 + Nh_K(P)).[K:\mathbb{Q}]^{(3 + Nh_K(P))}, $\newline

 $D:=\{4Nr.[K:\mathbb{Q}]^2.(H_K(P))^N.(N-1)!N! [\max_{i,j}H_{\bar{\mathbb{Q}}}(a_{i}^{(j)})]^2 \}^{2N^2r[K:\mathbb{Q}]^2}$,\newline \newline
it follows similarly as before that\newline \newline
$    \hat{h}_{\phi_{A}}(P) \geq \dfrac{1}{{\rho}^l l!}\sum_{v \in S} \sum_j(\sum_{i \leq r}{a}_{l_v}^{(i)}  \xi_id_{ij,v}(P))\log||x_j||_v +C_0\sum_{v \in R} \log||\zeta||_v~~~~~~~~ ~~~~~~~~~~~~~~~~~\newline
=\dfrac{1}{{\rho}^l l![K: \mathbb{Q}]}\sum_{v \in S} \sum_j(\sum_{i \leq r} D_{ij,v} \xi_i)\log||x_j||_v+ C_0\sum_{v \in R} \log||\zeta||_v ~~~~~~~~~~~~~~~~~~~~~~~~~~~~~~~~~~~~~~~~~~~~~~~~~~~~\newline
 =\dfrac{1}{{\rho}^l l!}\sum_{v \in S^0} [\sum_{i,j} D_{ij,v} \xi_i(-v(x_j))+\sum_{i,j,u \in S^{\infty}} D_{ij,u} \xi_i.v_{u,j}].\log N(v) + C_0\sum_{v \in R} \log||\zeta||_v,\newline \newline
> \dfrac{1}{2[K: \mathbb{Q}]}   \exp(-{E}.(\log A)^{(1+Nh_K(P))}.[\log D +  \log \log A]) + C_0m (\zeta, C),~~~~~~~~~~~~~~~~~~~~~~~~~~~~~~~~~~~~$\newline \newline

where $C= \dfrac{1}{2[K: \mathbb{Q}]}   \exp(-{E}.(\log A)^{(1+Nh_K(P))}.[\log D +  \log \log A]). $

\end{proof}

If $P \in \mathbb{G}_m^N(K)$ has algebraic integer coordinates, or in other words and notation, if $P \in (\mathcal{O}_K)^N$, then the constants $C_0$ and $h(\zeta)$ will not be needed for the bounds. \newline \newline
{\bf Theorem 3.4:}
\textit{Let $\phi:  {\mathbb{G}}_{m}^N \rightarrow  {\mathbb{G}}_{m}^N$ be a monomial map induced by $A$ with real eigenvalues, $l(A) \geq 1$ and $\delta_{\phi}>1$, and let $P \in  {\mathbb{G}}_{m}^N(\bar{\mathbb{Q}})$ be a point with orbit $\mathcal{O}_{\phi}(P)$ Zariski dense, such that the coordinates of $P$ in $\mathbb{Q}(J(A),P)$ are algebraic integers. Then there is an effective computable positive constant C depending on $\mathbb{Q}(J(A),P), A$ (more precisely on the heights of the coordinates of $J(A)$) and $h(P)$ such that $ \hat{h}_{\phi_{A}}(P)>C$}.
\begin{proof}
In this case $\zeta=1$ and thus $\log ||\zeta||_v=0$ for all $v$ on the proofs above.
\end{proof}

{\bf Corollary 3.5:}
\textit{Let $\phi:  {\mathbb{G}}_{m}^N \rightarrow  {\mathbb{G}}_{m}^N$ be a monomial map induced by $A$ with real eigenvalues, irreducible characteristic polynomial over $\mathbb{Q}$ and $\delta_{\phi}>1$, and let $P \in  {\mathbb{G}}_{m}^N(\bar{\mathbb{Q}})$ be a point with infinite orbit $\mathcal{O}_{\phi}(P)$, such that the coordinates of $P$ in $\mathbb{Q}(J(A),P)$ are algebraic integers. Then there is an effective computable positive constant C depending on $\mathbb{Q}(J(A),P), A$(more precisely on the heights of the coordinates of $J(A)$) and $h(P)$ such that $ \hat{h}_{\phi_{A}}(P)>C$}. \newline \newline
When the associated matrix is diagonalizable, Corollary 6.5 of [Lin, 11] shows that $\inf_n \dfrac{\deg(\phi^n)}{n^l \rho^n}\geq 1$, which yields the following case without $C_0$. \newline \newline
{\bf Corollary 3.6:}
\textit{Let $\phi:  {\mathbb{G}}_{m}^N \rightarrow  {\mathbb{G}}_{m}^N$ be a monomial map induced by $A$ diagonalizable with real eigenvalues, $l(A) \geq 1$ and $\delta_{\phi}>1$, and let $P \in  {\mathbb{G}}_{m}^N(\bar{\mathbb{Q}})$ be a point with orbit $\mathcal{O}_{\phi}(P)$ Zariski dense. Suppose $\zeta \in \mathbb{Z}$ satisfying that $\zeta P$ has algebraic integer coordinates is well-known. Then there is an effective computable positive constant C depending on $\mathbb{Q}(J(A),P), A$ (more precisely on the heights of the coordinates of $J(A)$), and on $h(P)$,  and  an effective function $m(\zeta, C )$ such that $ \hat{h}_{\phi_{A}}(P)>C+ m(\zeta, C )$}.  \newline \newline
In the general case of complex eigenvalues, the limit values $\xi_i$ in the proofs above aren't necessarily  $1$ or $-1$, and for instance we suppose they are algebraic numbers. If these are algebraic integers, since they have modulus 1, they are units, and thus have height 0, and the constants are the same as before.
\newline \newline
{\bf Theorem 3.7:}
\textit{Let $\phi:  {\mathbb{G}}_{m}^N \rightarrow  {\mathbb{G}}_{m}^N$ be a monomial map induced by $A$ with $l(A) \geq 1$ and $\delta_{\phi}>1$, and let $P \in  {\mathbb{G}}_{m}^N(\bar{\mathbb{Q}})$ be a point with orbit $\mathcal{O}_{\phi}(P)$ Zariski dense. Suppose $\zeta \in \mathbb{Z}$ satisfying that $\zeta P$ has algebraic integer coordinates, and $C_0=\inf_n \dfrac{\deg(\phi^n)}{n^l \rho^n} >0$, and also the limits numbers $\xi_i$'s, being algebraic numbers in the field $\mathbb{Q}(J(A),P)$(extended if necessary), are all well-known. Then there are an effective computable positive constant C depending on $\mathbb{Q}(J(A),P), A$ (more precisely on the heights of the coordinates of $J(A)$ and of the $\xi_i$'s), and $h(P)$, and  an effective function $m(\zeta , C)$ such that $ \hat{h}_{\phi_{A}}(P)>C+ C_0 m(\zeta, C )$}.
\begin{proof}
In a similar way as in the proof of Theorem 3.2, we have \newline \newline
$    \hat{h}_{\phi_{A}}(P) 
> \dfrac{1}{2{\rho}^l l![K: \mathbb{Q}]} \exp(-{E}^{\prime}.(\log A^{\prime})^{(1+Nh_K(P))}.[\log D^{\prime} +\log \log A^{\prime}])+ C_0m (\zeta, C),~~~~~~~~~~~~~~~~~~~~~~~~~~~~~~~~~~~~$ \newline \newline where \newline

$A^{\prime}::= e^{( Nh_K(P))(6.(1 + Nh_K(P)))}$,\newline

${E}^{\prime}:=C_{11}(1 + Nh_K(P)).[K:\mathbb{Q}]^{(3 + Nh_K(P))}, $\newline

 $D^{\prime}:= \{4Nr.[K:\mathbb{Q}]^2.(H_K(P))^N.(N-1)!N!\max_iH_{\bar{\mathbb{Q}}}(\xi_i)[\max_{i,j,l}H_{\bar{\mathbb{Q}}}(a_{il}^{(j)})]^2 \}^{2N^2r[K:\mathbb{Q}]^2}.$

\end{proof}

{\bf Corollary 3.8:}
\textit{Let $\phi:  {\mathbb{G}}_{m}^N \rightarrow  {\mathbb{G}}_{m}^N$ be a monomial map induced by $A$ with $l(A) \geq 1$ and $\delta_{\phi}>1$, and let $P \in  {\mathbb{G}}_{m}^N(\bar{\mathbb{Q}})$ be a point with orbit $\mathcal{O}_{\phi}(P)$ Zariski dense. Suppose $\zeta \in \mathbb{Z}$ satisfying that $\zeta P$ has algebraic integer coordinates, and $C_0=\inf_n \dfrac{\deg(\phi^n)}{n^l \rho^n} >0$, and also the limits numbers $\xi_i$'s, being algebraic integers in $\mathbb{Q}(J(A),P)$(extended if necessary), are all well-known. Then there are an effective computable positive constant C depending on $\mathbb{Q}(J(A),P), A$ (more precisely on the heights of the coordinates of $J(A)$), and $h(P)$, and  an effective function $m(\zeta , C)$ such that $ \hat{h}_{\phi_{A}}(P)>C+ C_0 m(\zeta, C )$}. \begin{proof}
As mentioned, in this case the $\xi_i$'s will have modulus $1$, and therefore will be units, with multiplicative height equal to $1$
\end{proof}
{\bf Remark 3.9:}  As Silverman did presume in [15, Remark 30], we see that the explicit constants obtained show that if the height of a point is very large, then the referred constants for this point will be smaller. In fact, expanding the expressions in the proof of Theorem 3.2, the constant  $C$ can be lower bounded by
\begin{center}
$ \dfrac{1}{2.\rho^l l![K: \mathbb{Q}] \{C_1.(H_K(P))^N.C(A)(1+h_K(P)))\}^{2N^2r.2^{53} \{2[K:\mathbb{Q}]^2.(1+ Nh_K(P))^2 \}^{(8(3 + N h_K(P)))}}} $,
\end{center}
where $C_1=24N^3(N-1)!^2r$ is an effective computable constant that does not depend on anything, except on $N$ eventually, and $C(A):=[\max_{i,j,l}H_{\bar{\mathbb{Q}}}(a_{il}^{(j)})]^2$ depends on the absolute heights of the coordinates of the vectors which transform $A$ in its Jordan form. As such vectors belong to the image of powers of $A$, the bigger the heights of coefficients of $A$ are, the smaller the constants obtained in the results are as well, roughly speaking. More sharp effective versions of effective Baker's theorem should lead to more sharp versions of the above constant.\newline \newline 
{\bf Remark 3.10:} Proposition 4 and Remark 5 of [12] show that bounds of Lehmer type for this canonical height do not exist in general.

\end{document}